\documentclass[11pt,draft]{article}
\usepackage{amssymb}
\usepackage{CJK}
\usepackage{amsmath,amsfonts,mathrsfs,amssymb}
\usepackage{indentfirst}

\numberwithin{equation}{section}

\begin{document}

\title{PMCV hypersurfaces in non-flat pseudo-Riemannian space forms}
\date{}
\maketitle
\vspace{-1.5cm}
\centerline{\author{Chao Yang$^1$, Jiancheng Liu$^{1,\ast}$ and Li Du$^2$}}
\begin{center}
\small{1. College of Mathematics and Statistics, Northwest Normal University, Lanzhou, 730070,
China\\
2. School of Science, Chongqing University of Technology,
Chongqing, 400054, China}
\end{center}

\vskip.3cm

\begin{quotation}
{\bf Abstract:}
{In this paper, we prove that PMCV (i.e. $\Delta\vec{H}$ is proportional to $\vec{H}$) hypersurface $M^n_r$ of a non-flat pseudo-Riemannian space form $N^{n+1}_s(c)$ with at most two distinct principal curvatures is minimal or locally isoparametric, and compute the mean curvature for the isoparametric ones. As an application, we give full classification results of such non-minimal Lorentzian hypersurfaces of non-flat Lorentz space forms.}



\end{quotation}

\footnote[0]{This work is supported by National Natural Science Foundation
of China (No. 12161078), Science and Technology Project of Gansu Province (No. 20JR5RA515), Innovation Ability Enhancement Project of Gansu Higher Education Institution (2023B-063).
}

\footnote[0]{E-mail address: liujc@nwnu.edu.cn}

\section{Introduction}

Let $N^{n+1}_s(c)$ be the ($n+1$)-dimensional pseudo-Riemannian space form
of constant curvature $c$ with index $s\geq 0$, and $x:\! M^n_r\to N^{n+1}_s(c)$ be an isometric immersion of
a pseudo-Riemannian hypersurface $M^n_r$ into $N^{n+1}_s(c)$.
Denote by $\vec{H}$ and $\Delta$ the mean curvature
vector field and the Laplace operator of $M^n_r$ with respect to the induced metric, respectively.

If the hypersurface $M^n_r$ satisfies the equation
$$
\Delta\vec{H}=\lambda\vec{H},
$$
for some real constant $\lambda$,
then it is said to have proper mean curvature vector field (abbreviated as PMCV).
Specially, when $\lambda=nc$, the hypersurface $M^n_r$ is known as biharmonic one.

The following geometric question was proposed by B.-Y. Chen in \cite{Chen 1988-1, Chen 2011}:

\emph{Determine all PMCV submanifolds of pseudo-Euclidean space $\mathbb{E}^{n+1}_s$. In particular, classify PMCV hypersurfaces in $\mathbb{E}^{n+1}_s$.}

In 1988, B.-Y. Chen \cite{Chen 1988-1} initiated the study of this problem for $s=0$, and proved that the PMCV surface $M^2$ in $\mathbb{E}^3$ is minimal, or locally congruent to a circular cylinder. In 1991, A. Ferr\'{a}ndez and P. Lucus showed in \cite{Ferrandez 1991} that the PMCV hypersurface $M^n$ of $\mathbb{E}^{n+1}$ with at most two distinct principal curvatures is minimal, or locally congruent to $\mathbb{E}^l\times\mathbb{S}^{n-l}(\frac{\lambda}{n-l})\ (0\leq l \leq n-1)$. Afterwards, B.-Y. Chen and O. J. Garay \cite{Chen-Garay 2012} derived that $\delta(2)$-ideal PMCV hypersurfaces in $\mathbb{E}^{n+1}$ is minimal, or locally congruent to $\mathbb{E}^1\times\mathbb{S}^{n-1}(\frac{\lambda}{n-1})$.

For the case $s=1$, A. Ferr\'{a}ndez and P. Lucus  \cite{Ferrandez 1992-1} classified PMCV surfaces, and generalized the result in \cite{Ferrandez 1992-2} to $n$-dimensional hypersurfaces whose minimal polynomial of the shape operator is at most of degree two. In \cite{Liu 2017}, C. Yang and J.-C. Liu weaken this condition to the number of distinct principal curvatures be not larger than two and obtained classification result.
For general index $s$, L. Du \cite{Du 2016} classified hypersuraces with diagonalizable shape operator and at most two distinct principal curvatures.

Naturally, it is interesting to consider the classification problem of PMCV hypersurfaces in non-flat pseudo-Riemannian space forms $N^{n+1}_s(c)$. Under the assumption that the shape operator is diagonalizable, L. Du \cite{Du 2016} classified such hypersurfaces of $N^{n+1}_s(c)\ (c\neq 0)$ with at most two distinct principal curvatures.
In this paper, we show that such PMCV hypersurface $M^n_r$ of $N^{n+1}_s(c)\ (c\neq 0)$ is minimal or locally isoparametric and give the value or range of the mean curvature, without the restriction that the shape operator is diagonalizable. As an application, we classify completely such non-minimal Lorentzian PMCV hypersurfaces of $N^{n+1}_1(c)\ (c\neq 0)$ by using classification theory of isoparametric hypersurfaces.

\section{Preliminaries}


Let $\mathbb{E}^{n}_s$\ $(0\leq s\leq n)$ be the $n$-dimensional pseudo-Euclidean space equipped with the canonical pseudo-Euclidean metric of index $s$
$$
g_0=-\sum_{i=1}^s\text{d}x_i^2+\sum_{j=s+1}^{n}\text{d}x_j^2,
$$
where $(x_1, x_2, \cdots, x_{n})$ is a rectangular coordinate system of $\mathbb{E}^{n}_s$.
Set
$$
\begin{aligned}
\mathbb{S}^n_s(c)=\left\{x\in\mathbb{E}^{n+1}_s\Big| -\sum_{i=1}^sx_i^2+\sum_{j=s+1}^{n+1}x_j^2=\frac{1}{c}>0\right\},\\
\mathbb{H}^n_s(c)=\left\{x\in\mathbb{E}^{n+1}_{s+1}\Big| -\sum_{i=1}^{s+1}x_i^2+\sum_{j=s+2}^{n+1}x_j^2=\frac{1}{c}<0\right\},
\end{aligned}
$$
where $c$ is a nonzero real constant. $\mathbb{S}^n_s(c)$ and $\mathbb{H}^n_s(c)$ are pseudo-Riemannian manifolds of curvature $c$ with index $s$, which together with $\mathbb{E}^{n}_s$ are called \emph{pseudo-Riemannian space forms}, denoted by $N^{n}_s(c)$. In particular, for $s=1$, $N^{n}_1(c)$ are called \emph{Lorentz space forms}, $\mathbb{S}^n_1(c)$ and $\mathbb{H}^n_1(c)$ known as \emph{de Sitter space} and \emph{anti-de Sitter space}, respectively.
A non-zero vector $X$ in $N^{n}_s(c)$ is said to be \emph{time-like},
\emph{space-like} or \emph{light-like}, according to whether $\langle X, X\rangle$ is negative, positive or zero.


Let $M^n_r$ be a nondegenerate hypersurface in $N^{n+1}_s(c)$,
$\vec{\xi}$ denote a unit normal vector field to $M^n_r$,
then $\varepsilon=\langle\vec\xi,\vec\xi\rangle=\pm1$.
Set $\nabla$ and $\tilde{\nabla}$ the Levi-Civita connections of $M^n_r$
and $N^{n+1}_s(c)$, respectively. For any vector fields $X, Y$ tangent to $M^n_r$,
it follows that
\begin{equation*}
\widetilde\nabla_{X}Y=\nabla_{X}Y + h(X,Y)\vec{\xi},
\end{equation*}
where $h$ is the scalar-valued second fundamental form.
Denote by $A$, $\vec{H}$ and $H$ the shape operator of $M^n_r$ associated to
 $\vec{\xi}$, the mean curvature vector field and the mean curvature,
 then $\vec{H}=H\vec{\xi}$ and $H=\frac{1}{n}\varepsilon\text{tr}A$.
For any vector fields $X, Y, Z$ tangent to $M^n_r$, the Codazzi and the Gauss equations are given by
\begin{equation}\label{Codazzi-eq}
\langle(\nabla_XA)Y, Z\rangle=
\langle(\nabla_YA)X, Z\rangle,
\end{equation}
\begin{equation*}
 R(X,Y)Z=c(\langle Y,Z\rangle X-\langle X,Z\rangle Y)+\varepsilon\langle A(Y),Z\rangle A(X)-\varepsilon\langle A(X),Z\rangle A(Y),
\end{equation*}
where $R(X,Y)Z=\nabla_{X}\nabla_{Y}Z-\nabla_{Y}\nabla_{X}Z-\nabla_{[X,Y]}Z$.

According to \cite{Du 2017}, $M^n_r$ is a PMCV hypersurface if and only if
\begin{equation}\label{proper-eq1}
 A(\nabla H)=-\frac{n}{2}\varepsilon H(\nabla H),
\end{equation}
\begin{equation}\label{proper-eq2}
\Delta H+\varepsilon H\text{tr}A^{2}=\lambda H.
\end{equation}

Suppose the eigenpolynomial $\text{det}(\lambda I-A)$ of $A$ can be expressed as $(\lambda-\lambda_1)^{n_1}(\lambda-\lambda_2)^{n_2}\cdots(\lambda-\lambda_k)^{n_k}$,
then $n_i$ is called the \emph{algebraic multiply} of $\lambda_i$. The \emph{geometric multiply} of $\lambda_i$ refers to the dimension of the eigensubspace
$\text{ker}(\lambda_i I-A):=\{X| AX=\lambda_iX\}$.
A hypersurface whose shape operator has constant eigenpolynomial
 is said to be \emph{isoparametric} (cf. \cite{Hahn 1984}),
which is equivalent to the principal curvatures and their algebriac multiplicities being constants.
In this paper, without special statement, multiplicity refers to algebraic multiplicity.

\section{PMCV hypersurfaces in non-flat pseudo-Riemannian space forms}

In this section, we study PMCV hypersurface $M^n_r$ of $N^{n+1}_s(c)\ (c\neq 0)$
with at most two distinct principal curvatures, and prove that it is minimal or locally isoparametric (see Theorem 3.1).
Furthermore, we estimate or compute the value of the mean curvature (see Theorems 3.3 and 3.4).

First of all, we recall from \cite{Chen 2011, Liu 2017} that the shape operator $A$ of the hypersurface $M^n_r$
can be expressed as an almost diagonal matrix
$$
A=\text{diag}\{A_{1}, \cdots, A_{t}, \bar{A}_{t+1}, \cdots, \bar{A}_{m}\},\ \text{for some}\ t,\ 0\leq t\leq m,
$$
where
\begin{small}
\begin{equation*}
\begin{aligned}
A_i\!=\!\left(
  \begin{array}{ccccc}
     \lambda_{i}&                &        &              &\\
               1&    \lambda_{i} &        &              &\\
                &          \ddots&  \ddots&              &\\
                &                &      1 &  \lambda_{i} &\\
                &                &        &            1 &      \lambda_{i}\\
  \end{array}
\right)\!,\
\bar{A}_j\!=\!\left(
  \begin{array}{ccccccc}
     \gamma_j&  \tau_j&          &         &           &              \\
     -\tau_j& \gamma_j&         0&         &           &         \\
      1&             0&  \gamma_j&  \tau_j &           &              \\
      &         \ddots&    \ddots&   \ddots&     \ddots&    \\
      &               &        1&          0& \gamma_j&    \tau_j\\
      &               &         &          1&  -\tau_j&   \gamma_j\\
  \end{array}
\right), \tau_j\neq 0,
\end{aligned}
\end{equation*}
\end{small}for $1\leq i\leq t$ and $t+1\leq j\leq m$, relative to the basis $\mathfrak{B}\!=\!\{u_{i_1}, u_{i_2}, \cdots, u_{i_{\alpha_i}}, u_{\bar{j}_1},$ $u_{\tilde{j}_1}, u_{\bar{j}_2}, u_{\tilde{j}_2}, \cdots, u_{\bar{j}_{\beta_j}},
u_{\tilde{j}_{\beta_j}}| 1\leq i\leq t,\ t+1\leq j\leq m\}$ with all scalar products zero except
$
\langle u_{i_a}, u_{i_b}\rangle=\varepsilon_{i}=\pm1$ and $\langle u_{\bar{j}_c},u_{\bar{j}_d}\rangle=1=-\langle u_{\tilde{j}_c}, u_{\tilde{j}_d}\rangle
$
if
$a+b=\alpha_i+1$. Here $1\leq i\leq t,\ c+d=\beta_j+1,\ t+1\leq j\leq m$ and
$
\sum_{i=1}^t\alpha_i+2\sum_{j=t+1}^m\beta_j=n.
$

It follows from the form of $A$ that $M^n_r$ has real principal curvatures
$
\lambda_1, \lambda_2, \cdots, \lambda_t
$
and imaginary eigenvalues
$
\gamma_{t+1}+\tau_{t+1}\sqrt{-1}, \gamma_{t+1}-\tau_{t+1}\sqrt{-1},
\cdots, \gamma_m+\tau_m\sqrt{-1}, \gamma_m-\tau_m\sqrt{-1}.
$
Notice that the imaginary principal curvatures appear as conjugate pairs.

\vskip.2cm
{\bf Theorem 3.1}\quad \emph{Let $M^n_r$ be a PMCV hypersurface of $N^{n+1}_s(c)\ (c\neq 0)$.
Suppose that $M^n_r$ has at most two distinct principal curvatures, then $M^n_r$ has constant mean curvature.
Furthermore, it is minimal or locally isoparametric.}

\vskip.2cm
{\bf Proof}\quad Firstly, we prove that if the mean curvature $H$ of $M^n_r$ is a constant,
then $M^n_r$ is minimal or locally isoparametric.
Let $x\in M^n_r$, there exists a neighborhood $U_x$ such that
the multiplicities of distinct principal curvatures are constants on $U_x$.
We work on $U_x$ and show that $M^n_r$ is isoparametric.

Suppose $\mu$ and $\nu$ are distinct principal curvatures of $M^n_r$ with multiplicities $l$ $(1\leq l\leq n)$ and $n-l$, then
\begin{equation}\label{trAlmun-lnu}
n\varepsilon H=\text{tr}A=l\mu+(n-l)\nu,
\end{equation}
and $\text{tr}A^2=l\mu^2+(n-l)\nu^2$. If $H$ is a nonzero constant, then \eqref{proper-eq2} implies $\text{tr}A^2=\varepsilon \lambda$. So,
\begin{equation}\label{varepsilonlambdalmun-lnu}
l\mu^2+(n-l)\nu^2=\varepsilon \lambda.
\end{equation}
It follows from \eqref{trAlmun-lnu} and \eqref{varepsilonlambdalmun-lnu} that $\mu$ and $\nu$ are constants. Hence, $M^n_r$ is isoparametric on $U_x$.

In the following, we show that $M^n_r$ has constant mean curvature.

When the principal curvatures are imaginary, $-\frac{n}{2}\varepsilon H$ (a real function) is not an eigenvalue of $A$.
It follows from \eqref{proper-eq1} that $\nabla H\equiv0$, i.e. $H=\text{const.}$
We need only treat the situation that the principal curvatures are all real.

Assume that $H$ is not a constant, then there exists a neighborhood $U_x$ of $x\in M^n_r$ such that $H\neq 0$ and $\nabla H\neq 0$.
From (\ref{proper-eq1}), we get that $-\frac{n}{2}\varepsilon H$ is an eigenvalue of $A$, and $\nabla H$ is a corresponding eigenvector.
When the principal curvatures of $M^n_r$ are all the same, denoted by $\mu$, then $\mu=\frac{1}{n}\text{tr}A=\varepsilon H$.
Thus, $-\frac{n}{2}\varepsilon H\ (\neq \mu)$ is not an eigenvalue of $A$, a contradiction.

Now, we assume that the number of distinct principal curvatures is two, i.e.
$A=\text{diag}\{A_{1}, \cdots, A_{m}\}$ and there are two distinct values among $\{\lambda_1, \cdots, \lambda_m\}$.
Observe the form of $A$, $\nabla H$ is in one of the directions $u_{i_{\alpha_i}}, 1\leq i\leq m$.
Without loss of generality, we suppose $\nabla H$ is in
the direction of $u_{1_{\alpha_1}}$, which may or may not be a light-like vector.
We follow different processes to lead contradictions for those two cases.

At the begining, we provide some important equations. We remark that, if no otherwise specified,
the ranges of $i$ and $a$ in the index $i_a$ are as:
$1\leq i\leq m$ and $1\leq a\leq \alpha_i$.
And when $a=\alpha_i$ (or $a=1$), the connection coefficients with the index $i_{a+1}$ (or $i_{a-1}$) disappear.

Let $[i]=\{j|\ \lambda_j=\lambda_i\}$ for $1\leq i\leq m$ and $\nabla_{u_{i_a}}u_{j_b}=
\sum_{k=1}^{m}\sum_{d=1}^{\alpha_k}\Gamma^{k_d}_{i_aj_b}u_{k_d}$. In view of the assumption that $\nabla H$ is in
the direction of $u_{1_{\alpha_1}}$, as well as the compatibility and symmetry of the connection, we obtain
\begin{equation}\label{B(H)-compatibility-eq12-symmetry-4case2}
\begin{cases}
u_{1_{1}}(H)\neq0,\ u_{B}(H)=0, \
\Gamma^{1_{1}}_{BC}=\Gamma^{1_{1}}_{CB},\\
\Gamma^{i_{\alpha_i-a+1}}_{Di_a}=0,\ \Gamma^{i_{\alpha_i-b+1}}_{Di_a}=-\Gamma^{i_{\alpha_i-a+1}}_{Di_{b}},\
\Gamma^{j_{\alpha_j-d+1}}_{Di_a}=
-\varepsilon_i\varepsilon_j\Gamma^{i_{\alpha_i-a+1}}_{Dj_{d}},
\end{cases}
\end{equation}
for $B, C, D\in\{k_e, 1\leq k\leq m, 1\leq e\leq\alpha_k\}$ and $B, C\neq 1_1$.

Putting $X=u_{i_a}, Y=u_{j_b}, Z=u_{k_{\alpha_k-d+1}}$ into \eqref{Codazzi-eq}, it gives
\begin{equation}\label{Codazzi-td-eq}
\begin{aligned}
u_{i_a}\!(\lambda_j)\delta_{j_b}^{k_d}\!+\!(\lambda_j-&\lambda_k)\Gamma_{i_aj_b}^{k_d}\!+\!\Gamma_{i_aj_{b+1}}^{k_d}\!-\!\Gamma_{i_aj_b}^{k_{d-1}}
\!=\!u_{j_b}\!(\lambda_i)\delta_{i_a}^{k_d}\!+\!(\lambda_i-\lambda_k)\Gamma_{j_bi_a}^{k_d}\!+\!\Gamma_{j_bi_{a+1}}^{k_d}\!-\!\Gamma_{j_bi_a}^{k_{d-1}},
\end{aligned}
\end{equation}
where $\delta_{j_b}^{i_a}=1$ if $i_a=j_b$, and $\delta_{j_b}^{i_a}=0$ if $i_a\neq j_b$.
Using \eqref{B(H)-compatibility-eq12-symmetry-4case2}, we deduce from \eqref{Codazzi-td-eq} that
\begin{equation}\label{codazzi-eqs}
\begin{aligned}
&\text{for}\ i=j=k=d=1\Rightarrow \Gamma_{1_{\alpha_1-1}1_{\alpha_1}}^{1_{\alpha_1-1}}\!=\cdots=\!\Gamma_{1_11_{\alpha_1}}^{1_1}\!=0;\\
&\text{for}\ i=k=d=1, j\in[1], j\neq 1\Rightarrow \Gamma_{1_{a+1}j_{\alpha_j}}^{1_1}\!=0;\\
&\text{for}\ i=k=d=1, j\notin[1]\Rightarrow \Gamma_{1_{a}j_b}^{1_1}\!=0;\\
&\text{for}\ i\notin[1], j=1, b=\alpha_1, k=1\Rightarrow \Gamma_{1_{\alpha_1}i_a}^{1_d}\!=0;\\
&\text{for}\ i,j\neq 1, j \notin [i], k=d=1 \Rightarrow  \Gamma_{i_aj_b}^{1_1}\!=0;\\
&\text{for}\ i,j\!\neq\! 1, j\!\in\![i], k=\!d\!=\!1 \!\Rightarrow\!\Gamma_{i_a1_{\alpha_1}}^{j_b} \!=\!\Gamma_{i_{a+1}1_{\alpha_1}}^{j_{b+1}}\!,\ \Gamma_{i_{\alpha_i}j_{b+1}}^{1_1}\!=\!\Gamma_{i_a\!1_{\alpha_1}\!}^{j_e}\!=\!0, e<a;\\
&\text{for}\ \!i\!=\!k(\!\neq \!1)\!, d\!=\!a, j\!=\!1, b\!=\!\alpha_1\!\Rightarrow\! u_{1_{\alpha_1}}\!(\lambda_i)\!=\!(-\frac n2\varepsilon H-\lambda_i)\Gamma_{i_{\alpha_i}1_{\alpha_1}}^{i_{\alpha_i}}\!,\ \Gamma_{1_{\alpha_1}i_{a+1}}^{i_a}\!=\!0;\!\\
&\text{for}\ i,a\neq 1, j=d=1, b=\alpha_1, k=i\Rightarrow \Gamma_{1_{\alpha_1}i_{a}}^{i_1}\!=0.\\
\end{aligned}
\end{equation}

Let
\begin{equation*}
\begin{cases}
A_{i,0,j_b}=\Gamma_{i_1j_b}^{i_1}
+\Gamma_{i_2j_b}^{i_2}+\cdots
+\Gamma_{i_{\alpha_i}j_b}^{i_{\alpha_i}},\\
A_{i,1,j_b}=\Gamma_{i_2j_b}^{i_1}
+\Gamma_{i_3j_b}^{i_2}+\cdots
+\Gamma_{i_{\alpha_i}j_b}^{i_{\alpha_i-1}},\\
\quad\quad\quad\quad\quad\vdots\\
A_{i,\alpha_i-1,j_b}=\Gamma_{i_{\alpha_i}j_b}^{i_1},
\end{cases}\ \ \text{for} \ j\in[1],\ i\notin[1],
\end{equation*}
then \eqref{Codazzi-td-eq} for $k=i$ implies
\begin{equation*}
\begin{cases}
(-\frac n2\varepsilon H-\lambda_i)A_{i,0,j_b}
+A_{i,0,j_{b+1}}-A_{i,1,j_b}=\alpha_iu_{j_b}(\lambda_i),\\
(-\frac n2\varepsilon H-\lambda_i)A_{i,r,j_b}+A_{i,r,j_{b+1}}-A_{i,r+1,j_b}=0,\\
(-\frac n2\varepsilon H-\lambda_i)A_{i,\alpha_i-1,j_b}+A_{i,\alpha_i-1,j_{b+1}}=0,
\end{cases}
\end{equation*}
$1\leq r\leq \alpha_i-2$. It follows that
\begin{equation}\label{A0case2-----1}
A_{i,s,j_b}=0,\ 1\leq s\leq \alpha_i-1,\\
\end{equation}
\begin{equation}\label{Ano0case2}
(-\frac n2\varepsilon H-\lambda_i)A_{i,0,j_b}
+A_{i,0,j_{b+1}}=\alpha_iu_{j_b}(\lambda_i).
\end{equation}
Set
\begin{equation*}
\begin{cases}
B_{i,k,0,j_b}=\Gamma_{i_1j_b}^{k_1}
+\Gamma_{i_2j_b}^{k_2}+\cdots
+\Gamma_{i_{y_0}j_b}^{k_{y_0}},\\
B_{i,k,1,j_b}=\Gamma_{i_2j_b}^{k_1}
+\Gamma_{i_3j_b}^{k_2}+\cdots
+\Gamma_{i_{y_1}j_b}^{k_{y_1-1}},\\
\quad\quad\quad\quad\quad\vdots\\
B_{i,k,\alpha_i-1,j_b}=\Gamma_{i_{\alpha_i}j_b}^{k_1},
\end{cases}\!
\end{equation*}
for $j\in[1]$, $i, k\notin[1]$, $i\neq k$ and $y_r=\min \{\alpha_i, \alpha_k+r\}$, $0\leq r\leq \alpha_i-1$, then we obtain from \eqref{Codazzi-td-eq} that
\begin{equation*}
\begin{cases}
(-\frac n2\varepsilon H-\lambda_k)B_{i,k,s,j_b}+B_{i,k,s,j_{b+1}}
-B_{i,k,s+1,j_b}=0,\\
(-\frac n2\varepsilon H-\lambda_k)B_{i,k,\alpha_i-1,j_b}+B_{i,k,\alpha_i-1,j_{b+1}}=0,
\end{cases}
\end{equation*}
which implies
\begin{equation}\label{A0case2-----2}
B_{i,k,r,j_b}=0.
\end{equation}

\vskip.3cm
We treat the two cases that $\nabla H$ is not light-like or light-like respectively and deduce contradictions.

\vskip.3cm
\emph{Case 1:\ $\nabla H$ is not light-like.}



Since $\nabla H$ is in the direction $u_{1_{\alpha_1}}$, we have $\langle u_{1_{\alpha_1}}, u_{1_{\alpha_1}}  \rangle\neq 0$, which implies $2\alpha_1=\alpha_1+1$, i.e. $\alpha_1=1$.
We find from the seventh equation of \eqref{codazzi-eqs} that
if $\lambda_i=-\frac n2\varepsilon H$ for some $i\neq 1$, then $u_{1_1}(H)=0$,
which contradicts to \eqref{B(H)-compatibility-eq12-symmetry-4case2}.
Hence, $\lambda_i\neq -\frac n2\varepsilon H$ for $2\leq i\leq m$.
Combining $\lambda_1=-\frac n2\varepsilon H$ and $\text{tr}A=n\varepsilon H$, we know $\lambda_2=\cdots=\lambda_m=\frac{3n\varepsilon H}{2(n-1)}$.
It follows from the seventh equation of \eqref{codazzi-eqs} that
\begin{equation}\label{W--3.1}
W:=\Gamma^{i_{a}}_{i_{a}1_{1}}\!=-\frac{3u_{1_1}(H)}{(n+2)H},\ 2\leq i\leq m.
\end{equation}

\vskip.2cm\
{\bf Lemma 3.2}\quad \emph{We have}
\begin{equation}\label{i--3.1case1}
u_{1_1}(W)+W^2=\frac{3\varepsilon\varepsilon_1 n^2H^2}{4(n-1)}-c\varepsilon_1.
\end{equation}

{\bf Proof}\quad 
Let $j=b=1, d=1$ in \eqref{Codazzi-td-eq}, then
\begin{equation}\label{bu66-}
-\frac{n(n+2)\varepsilon H}{2(n-1)}\Gamma_{i_{a}1_1}^{k_1}\!=
\Gamma_{1_1i_{a+1}}^{k_1},\ 2\leq i, k\leq m,\ i\neq k.
\end{equation}
Combining the sixth equation of \eqref{codazzi-eqs}, it leads to
\begin{equation}\label{57}
\Gamma_{1_1i_{a+1}}^{k_1}\!=0,\ 2\leq a\leq\alpha_i-1.
\end{equation}

Calculate $\langle R(u_{1_1}, u_{i_{\alpha_i}})u_{1_1},u_{i_1}\rangle$ by using the Gauss equation, combining \eqref{codazzi-eqs} and \eqref{57}, we have
\begin{equation}\label{PW}
u_{1_1}(W)+W^2-\frac{3\varepsilon\varepsilon_1 n^2H^2}{4(n-1)}+c\varepsilon_1=P_i,\ 2\leq i\leq m,
\end{equation}
where
$$
P_i=
\sum_{2\leq k\leq m, k\neq i}\!\left(\Gamma_{k_{\alpha_k}1_1}^{i_{\alpha_i}}
\Gamma_{1_1i_{\alpha_i}}^{k_{\alpha_k}}\!-\Gamma_{i_{\alpha_i}1_1}^{k_{\alpha_k}}
\Gamma_{1_1k_{\alpha_k}}^{i_{\alpha_i}}\!
+\Gamma_{k_{\alpha_k-1}1_1}^{i_{\alpha_i}}
\Gamma_{1_1i_{\alpha_i}}^{k_{\alpha_k-1}}\!
-\Gamma_{i_{\alpha_i}1_1}^{k_{\alpha_k}}
\Gamma_{k_{\alpha_k}1_1}^{i_{\alpha_i}}\right).
$$

From the sixth equation of \eqref{codazzi-eqs}, we obtain that $\Gamma^{i_{1}}_{k_{1}1_1}\!=\Gamma^{i_{\alpha_k}}_{k_{\alpha_k}1_1}\!=0$ if $\alpha_k<\alpha_i\ (i, k\neq 1)$.
And if $\alpha_k=\alpha_i\ (i, k\neq 1)$,
\eqref{A0case2-----2} tells us that
$
\Gamma^{i_{1}}_{k_{1}1_1}\!+\Gamma^{i_{2}}_{k_{2}1_1}\!+\cdots+\Gamma^{i_{\alpha_i}}_{k_{\alpha_k}1_1}\!=0,
$
which together with the sixth equation of \eqref{codazzi-eqs} implies $\Gamma^{i_{1}}_{k_{1}1_1}\!=0$.
If $\alpha_k>\alpha_i\ (i, k\neq 1)$, taking $j=b=1,\ d=a+1$ in \eqref{Codazzi-td-eq} and combining
$\Gamma^{k_{1}}_{i_{1}1_1}\!=0$, we have
\begin{equation}\label{chzh--2---lm3.6}
-\frac{n(n+2)\varepsilon H}{2(n-1)}\Gamma^{k_{a+1}}_{i_{a}1_1}=
\Gamma^{k_{a+1}}_{1_1i_{a+1}}-\Gamma^{k_{a}}_{1_1i_{a}}.
\end{equation}
When $\alpha_k>\alpha_i+1$, combining \eqref{B(H)-compatibility-eq12-symmetry-4case2}, the sixth equation of \eqref{codazzi-eqs} implies that
\begin{equation}\label{i11alpha1j2=0---lm3.6}
\Gamma^{k_{2}}_{i_{1}1_1}\!=\cdots=\Gamma^{k_{\alpha_i+1}}_{i_{\alpha_i}1_1}\!=0.
\end{equation}
\eqref{chzh--2---lm3.6} together with \eqref{i11alpha1j2=0---lm3.6} yields
\begin{equation}\label{1alpha1i1j1--2---lm3.6}
\Gamma^{k_{1}}_{1_1i_{1}}\!=0,\ \text{if}\ \alpha_k>\alpha_i+1.
\end{equation}
When $\alpha_k=\alpha_i+1$, we deduce from \eqref{B(H)-compatibility-eq12-symmetry-4case2}, the sixth equation of \eqref{codazzi-eqs} and
\eqref{chzh--2---lm3.6} that
\begin{equation}\label{buzhi--1---lm3.6}
\begin{aligned}
\Gamma^{k_{2}}_{i_{1}1_1}\!=\!\Gamma^{k_{\alpha_k}}_{i_{\alpha_i}1_1}
\!=\!\varepsilon_i\varepsilon_k\Gamma^{i_{1}}_{k_{1}1_1}\!,\
\Gamma^{k_{1}}_{1_1i_{1}}\!=\!\alpha_i\Gamma^{k_{\alpha_i}}_{1_1i_{\alpha_i}}
\!=\!\varepsilon_i\varepsilon_k\frac{n(n+2)\varepsilon H}{2(n-1)}\alpha_i\Gamma^{i_{1}}_{k_{1}1_1}\!.
\end{aligned}
\end{equation}

By \eqref{B(H)-compatibility-eq12-symmetry-4case2},
we know $\Gamma_{i_{\alpha_i}1_1}^{j_{\alpha_j}}=\varepsilon_i\varepsilon_j\Gamma_{j_{1}1_1}^{i_{1}}$
and $\Gamma_{j_{\alpha_j-1}1_1}^{i_{\alpha_i}}=\varepsilon_i\varepsilon_j\Gamma_{i_{1}1_1}^{j_{2}}$.
Using the sixth equation of \eqref{codazzi-eqs}, \eqref{bu66-}, \eqref{i11alpha1j2=0---lm3.6}, \eqref{1alpha1i1j1--2---lm3.6} and \eqref{buzhi--1---lm3.6},
$P_i$ can be rewritten as
\begin{equation}\label{i--2----lm3.6}
\begin{aligned}
P_i=&
\!\sum_{k\neq i}\!\left(-\Gamma_{i_{1}1_1}^{k_{1}}
\Gamma_{1_1k_{1}}^{i_{1}}\!-\Gamma_{i_{1}1_1}^{k_{2}}
\Gamma_{1_1k_{2}}^{i_{1}}\!
+\Gamma_{k_{1}1_1}^{i_{1}}
\Gamma_{1_1i_{1}}^{k_{1}}\!
-\Gamma_{i_{1}1_1}^{k_{1}}
\Gamma_{k_{1}1_1}^{i_{1}}\right)\\
=&
\!-\!\sum_{\alpha_k=\alpha_i-1}\!\Gamma_{i_{1}1_1}^{k_{1}}
\Gamma_{1_1k_{1}}^{i_{1}}\!-\!\sum_{\alpha_k=\alpha_i+1}\!\Gamma_{i_{1}1_1}^{k_{2}}
\Gamma_{1_1k_{2}}^{i_{1}}\!
+\!\sum_{\alpha_k=\alpha_i+1}\!\Gamma_{k_{1}1_1}^{i_{1}}
\Gamma_{1_1i_{1}}^{k_{1}}\\
=&
\frac{n(n+2)\varepsilon H}{2(n-1)}\varepsilon_i\left[\sum_{\alpha_k= \alpha_i+1}\!\varepsilon_k(\alpha_i+1)(\Gamma^{i_{1}}_{k_{1}1_{1}})^2-\!\sum_{\alpha_k= \alpha_i-1}\!\varepsilon_k\alpha_k(\Gamma^{k_{1}}_{i_{1}1_{1}})^2\right].
\end{aligned}
\end{equation}
Choose $\{\alpha_{1_1},\alpha_{1_2},\cdots,\alpha_{1_{\xi_1}},\alpha_{2_1},\alpha_{2_2},\cdots,\alpha_{2_{\xi_2}},\cdots,\alpha_{h_1},\alpha_{h_2},\cdots,\alpha_{h_{\xi_h}}\}$ in $\{\alpha_2,\cdots,\alpha_m\}$\\
such that $
\tilde{\alpha}_t:=\alpha_{t_1}=\alpha_{t_2}=\cdots=\alpha_{t_{\xi_t}}$ for $1\leq t\leq h$,\   $\tilde{\alpha}_k+1=\tilde{\alpha}_{k+1}$ for $1\leq k\leq h-1$,
and $\tilde{\alpha}_{1}-1, \tilde{\alpha}_h+1\notin\{\alpha_2,\cdots,\alpha_m\}$.
Since $P_i=P_k$ for distinct $i, k$, we can denote $P_i$ by $P$. Puting $i=t_1,t_2, \cdots, t_{\xi_t}$ in \eqref{i--2----lm3.6} and taking the sum of those equations gives
\begin{equation}\label{i--3----lm3.6}
\xi_tP=(\tilde{\alpha}_{t}+1)Q_{t}-\tilde{\alpha}_{t-1}Q_{t-1},\ 1\leq t\leq h,
\end{equation}
where
$$
Q_t=\!\sum_{{\!\tiny\!\begin{array}{c}\! 1\!\leq\! p\!\leq \!\xi_t\!\\ \!1\!\leq\! q\!\leq \!\xi_{t+1}\!\end{array}\!}\!}\!\frac{n(n+2)\varepsilon H}{2(n-1)}\varepsilon_{\rho}\varepsilon_{\sigma}(\Gamma^{\rho_{1}}_{\sigma_{1}1_{\alpha_1}})^2,\ \ \rho=t_p,\ \sigma=(t+1)_q.
$$
\eqref{i--3----lm3.6} implies $Q_t=0$, $1\leq t\leq h$
and $P=0$. Thus, \eqref{i--3.1case1} follows from \eqref{PW}.
\hfill$\square$

Now, we continue the proof of Theorem 3.1 for case 1.

Following the process of the proof of \cite[Lemma 3.4]{Liu 2017}, we derive
\begin{equation}\label{proper-td1-3.2.1}
\begin{aligned}
u_{1_1}u_{1_1}(H)+(n-1)Wu_{1_1}(H)
-\varepsilon_1\varepsilon \frac{(n+8)n^2H^3}{4(n-1)}+\lambda\varepsilon_1H=0.
\end{aligned}
\end{equation}
We obtain from \eqref{W--3.1} and \eqref{i--3.1case1} that
\begin{equation*}
\begin{aligned}
u_{1_1}u_{1_1}(H)=\frac{(n+2)(n+5)H}{9}W^2
-\frac{(n+2)n^2\varepsilon\varepsilon_1H^3}{4(n-1)}+\frac{1}{3}(n+2)c\varepsilon_1H.
\end{aligned}
\end{equation*}
Substituting into \eqref{proper-td1-3.2.1} leads to
\begin{equation}\label{proper-td3--3.2.1}
\begin{aligned}
\frac{2(n-4)(n+2)}{9}W^2
+\frac{n^2(n+5)\varepsilon\varepsilon_1H^2}{2(n-1)}-\frac{1}{3}(n+2)c\varepsilon_1-\varepsilon_1\lambda=0.
\end{aligned}
\end{equation}
Differentiating (\ref{proper-td3--3.2.1}) along $u_{1_1}$, and applying \eqref{W--3.1} and \eqref{i--3.1case1},
we deduce
\begin{equation*}
\begin{aligned}
-4(n-4)W^2+\frac{3 n^2}{n-1}\varepsilon\varepsilon_1H^2-4(n-4)c\varepsilon_1=0,
\end{aligned}
\end{equation*}
which together with (\ref{proper-td3--3.2.1}) gives a polynomial equation of $H$, which tells us that $H$ must be a constant, a contradiction.

\vskip.2cm
\emph{Case 2:\ $\nabla H$ is light-like.}

We know $\alpha_1\geq 2$ in this case.
Suppose $\lambda_1=\cdots=\lambda_p= -\frac{n}{2}\varepsilon H$
and $\lambda_{p+1}=\cdots=\lambda_m(\neq\!-\frac{n}{2}\varepsilon H)$. Denote $l=\alpha_1+\cdots+\alpha_p$,
then $\alpha_{p+1}+\cdots+\alpha_m=n-l$ since $\alpha_1+\cdots+\alpha_m=n$.
Combining $\sum_{i=1}^m\alpha_i\lambda_i=\text{tr}A=n\varepsilon H$, we get
$\lambda_m=\frac{(2+l)n\varepsilon H}{2(n-l)}$.
It follows from the seventh equation of \eqref{codazzi-eqs} and \eqref{Ano0case2} that
\begin{equation}\label{ukdtrA}
\begin{aligned}
\Gamma_{i_11_{\alpha_1}}^{i_1}\!=0,\ \ \ \sum_{a=1}^{\alpha_i}\left(\frac{-n(n+2)\varepsilon H}{2(n-l)}
\Gamma_{i_ak_d}^{i_a}
+\Gamma_{i_ak_{d+1}}^{i_a}\right)=0,
\end{aligned}
\end{equation}
for $2\leq i\leq m$, $1\leq k\leq p$ and $k_d\neq 1_1$. As $\Gamma_{i_a1_{\alpha_1}}^{j_b}\!=\Gamma_{i_{a+1}1_{\alpha_1}}^{j_{b+1}}$ for $j\in[i]$, $i, j\neq 1$ (see \eqref{codazzi-eqs}), we derive from \eqref{A0case2-----1} and \eqref{A0case2-----2} that
\begin{equation*}\label{3.21}
\begin{aligned}
\sum_{a=2}^{\alpha_i}\sum_{
1\leq b< a}
\Gamma_{i_a 1_{d}}^{i_b}\Gamma_{i_b1_{\alpha_1}}^{i_a}\!=\Gamma_{i_11_{\alpha_1}}^{i_2}A_{i,1,1_d}+\Gamma_{i_11_{\alpha_1}}^{i_3}A_{i,2,1_d}+\cdots+\Gamma_{i_11_{\alpha_1}}^{i_{\alpha_i}}A_{i,\alpha_i-1,1_d}=0,\\
\sum_{a=1}^{\alpha_i}\!\sum_{\tiny\begin{array}{c}
\!1\!\leq\! b\!\leq\!\alpha_k\!\\\!
b\leq a\!\end{array}}
\!\Gamma_{i_a 1_{d}}^{k_b}\Gamma_{k_b1_{\alpha_1}}^{i_a}\!=\!\Gamma_{k_11_{\alpha_1}}^{i_1}\!B_{i,k,0,1_d}\!+\!\Gamma_{k_11_{\alpha_1}}^{i_2}\!B_{i,k,1,1_d}\!+\!\cdots+\!\Gamma_{k_11_{\alpha_1}}^{i_{\alpha_i}}\!B_{i,k,\alpha_i-1,1_d}\!=\!0,
\end{aligned}
\end{equation*}
for $p+1\leq i, k\leq m$ and $i\neq k$.

Using the Gauss equation for $\langle R(u_{1_{2}}, u_{i_{a}})u_{1_{\alpha_1}}\!,u_{i_{\alpha_i-a+1}}\rangle$ and
$\langle R(u_{1_{1}}, u_{i_{a}})u_{1_{\alpha_1}}\!,u_{i_{\alpha_i-a+1}}\rangle$ for $p+1\leq i\leq m$,
and taking sum of the results over $i$ and $a$, as well as combining \eqref{B(H)-compatibility-eq12-symmetry-4case2}, \eqref{codazzi-eqs}, \eqref{ukdtrA} and the above two equations,
we deduce
\begin{equation}\label{a1----Lemma 3.14}
\begin{aligned}
\sum_{i=p+1}^m\left(\sum_{d=3}^{\alpha_1-1}
\Gamma_{1_21_{\alpha_1}}^{1_d}\!\sum_{a=1}^{\alpha_i}\!\Gamma_{i_a1_d}^{i_a}
+\sum_{k=2}^p\sum_{d=2}^{\alpha_k}\Gamma_{1_21_{\alpha_1}}^{k_d}
\!\sum_{a=1}^{\alpha_i}\!\Gamma_{i_ak_d}^{i_a}\right)=0,
\end{aligned}
\end{equation}
\begin{equation}\label{a2----Lemma 3.14}
\begin{aligned}
\sum_{i=p+1}^m\!\left(\sum_{d=2}^{\alpha_1-1}
\!\Gamma_{1_11_{\alpha_1}}^{1_d}\!\sum_{a=1}^{\alpha_i}\!\Gamma_{i_a1_d}^{i_a}
\!+\!\sum_{k=2}^p\sum_{d=1}^{\alpha_k}\!\Gamma_{1_11_{\alpha_1}}^{k_d}\!\sum_{a=1}^{\alpha_i}\!
\Gamma_{i_ak_d}^{i_a}\right)\!=\!
(n-l)c\varepsilon_1\!-
\frac{1}{4}(l+2) n^2\varepsilon\varepsilon_1H^2.
\end{aligned}
\end{equation}
Multiplying both sides of \eqref{a2----Lemma 3.14} by $\frac{n(n+2)\varepsilon H}{2(n-l)}$,
and combining \eqref{ukdtrA}, it gives
\begin{equation}\label{kl}
\begin{aligned}
\sum_{i=p+1}^m\left(\sum_{d=2}^{\alpha_1-2}\Gamma_{1_11_{\alpha_1}}^{1_d}
\sum_{a=1}^{\alpha_i}
\Gamma_{i_a1_{d+1}}^{i_a}
+\sum_{k=2}^p\sum_{d=1}^{\alpha_k-1}\Gamma_{1_11_{\alpha_1}}^{k_d}
\sum_{a=1}^{\alpha_i}\Gamma_{i_ak_{d+1}}^{i_a} \right)\\
=\frac{n(n+2)\varepsilon H}{2}\left(
c\varepsilon_1-\frac{(l+2)n^2\varepsilon\varepsilon_1H^2}{4(n-l)}\right).
\end{aligned}
\end{equation}
By \eqref{a1----Lemma 3.14},
equation \eqref{kl} reduces to
$$
4(n-l)c\varepsilon_1-(2+l)n^2\varepsilon\varepsilon_1H^2=0.
$$
Therefore, $H$ is a constant, a contradiction.
$\hfill\square$

\vskip.2cm
Based on the results of Theorem 3.1, in the following, we estimate the range of that constant mean curvature (see Theorem 3.3),
and compute that constant for the special case that the the algebraic and geometric multiplicities of some principal curvature coincide (see Theorem 3.4).

We illustrate that under the assumption that the number of distinct principal curvatures is not larger than two, the principal curvatures are all real or imaginary, since imaginary principal curvatures are appear as conjugate pairs. As before, $\varepsilon$ is the inner product of the normal vector field $\xi$ with itself.

\vskip.2cm
{\bf Theorem 3.3}\quad \emph{Let $M^n_r$ be a non-minimal PMCV hypersurface of $N^{n+1}_s(c) (c\neq 0)$, with PMCV constant $\lambda$.
\begin{itemize}
  \item When $M^n_r$ has one or two distinct real principal curvatures, we have $
H^2\leq\frac{\varepsilon\lambda}{n}$ with $\varepsilon \lambda>0$. The equality holds if and only if the principal curvatures are all equal;
  \item When $M^n_r$ has a pair of imaginary principal curvatures, we have $H^2>\frac{\varepsilon\lambda}{n}$.
\end{itemize}
}

{\bf Proof}\quad Suppose $\mu$ and $\nu$ are distinct principal curvatures of $M^n_r$, with multiplicities $l$ $(1\leq l\leq n)$ and $n-l$.

When $\mu$ and $\nu$ are real, we get from \eqref{trAlmun-lnu} and \eqref{varepsilonlambdalmun-lnu} that $
n^2H^2=(l\mu+(n-l)\nu)^2\leq n(l\mu^2+(n-l)\nu^2)=n\varepsilon \lambda
$. So, we have $\varepsilon \lambda>0$ and $H^2\leq\frac{\varepsilon\lambda}{n}$, with the equality holding if and only if $l=n$.

When $\mu$ and $\nu$ are a pair of imaginary numbers with multiply $l=\frac{n}{2}$, we suppose $\mu=\gamma+\tau\sqrt{-1}$ and $\nu=\gamma-\tau\sqrt{-1}$, with $\tau\neq 0$. Then \eqref{trAlmun-lnu} and \eqref{varepsilonlambdalmun-lnu} reduces to $\gamma=\varepsilon H$ and $n(\gamma^2-\tau^2)=\varepsilon \lambda$. It follows that $n H^2>n H^2-n\tau^2=n(\gamma^2-\tau^2)=\varepsilon \lambda$, i.e. $H^2>\frac{\varepsilon \lambda}{n}$.
$\hfill\square$

\vskip.2cm
{\bf Remark 3.4}\quad We point out that the case that $M^n_r$ has a pair of imaginary principal curvatures in Theorem 3.3 does not exist when $n$ is an odd number.

\vskip.2cm
In particular, when the algebraic and geometric multiplicities of some principal curvature coincide, we can compute the value of $H$.



\vskip.2cm
{\bf Theorem 3.5}\quad \emph{Let $M^n_r$ be a non-minimal PMCV hypersurface of $N^{n+1}_s(c) (c\neq 0)$, with PMCV constant $\lambda$. Suppose that $M^n_r$ has two distinct principal curvatures $\mu$ (multiply $l$) and $\nu$, and the algebraic and geometric multiplicities of $\mu$ or $\nu$ coincide, then one of the following holds.
\begin{itemize}
  \item  When $\mu$ and $\nu$ are real, we have
\begin{equation*}
H^2=\frac{1}{2n^2}\left[n\varepsilon\lambda-4l(n-l)c\varepsilon\pm(2l-n)\sqrt{\lambda^2-4l(n-l)c^2}\right],
\end{equation*}
$$
\mu^2=\frac{\varepsilon\lambda\pm \sqrt{\lambda^2-4l(n-l)c^2}}{2l},\ \nu^2=\frac{\varepsilon\lambda\mp\sqrt{\lambda^2-4l(n-l)c^2}}{2(n-l)},
$$
with $\varepsilon \lambda\geq2\sqrt{l(n-l)}|c|$. In particular, if $\varepsilon\lambda=n|c|$, then $H^2=-\frac{ 4l(n-l)c\varepsilon}{n^2}
$, $\mu^2=-\frac{(n-l) c\varepsilon}{l}$, $\nu^2=-\frac{l c\varepsilon}{n-l}
$ for $c\varepsilon<0$, or $H^2=\frac{ (2l-n)^2c\varepsilon}{n^2}
$, $\mu^2=\nu^2=c\varepsilon $ for $c\varepsilon>0$.
    \item When $\mu$ and $\nu$ are imaginary, we have $H^2=\frac{1}{2n}(\varepsilon\lambda-nc\varepsilon)$ with $c\varepsilon<0$ and $|\varepsilon\lambda|<-nc\varepsilon$.
\end{itemize}
}


{\bf Proof}\quad According to Theorem 3.1, $M^n_r$ is locally isoparametric. Then, by the basic identity of Cartan in \cite[Theorem 2.9]{Hahn 1984}, we obtain $c+\varepsilon\mu\nu=0$.

(i)\ \ When $\mu$ and $\nu$ are real, by applying \eqref{trAlmun-lnu}, \eqref{varepsilonlambdalmun-lnu} and $c+\varepsilon\mu\nu=0$, we have
\begin{equation*}
\begin{cases}
\varepsilon \lambda+2\sqrt{l(n-l)}c\varepsilon=(\sqrt{l}\mu-\sqrt{n-l}\nu)^2\geq 0,\\
\varepsilon \lambda-2\sqrt{l(n-l)}c\varepsilon=(\sqrt{l}\mu+\sqrt{n-l}\nu)^2\geq 0,
\end{cases}
\end{equation*}
which implies $\varepsilon\lambda\geq2\sqrt{l(n-l)}|c|$.
Calculate $\mu^2$, $\nu^2$ and $H^2$ from \eqref{trAlmun-lnu}, \eqref{varepsilonlambdalmun-lnu} and $c+\varepsilon\mu\nu=0$, we derive
\begin{equation*}\label{mu2nu2}
\mu^2=\frac{\varepsilon\lambda\pm\sqrt{\lambda^2-4l(n-l)c^2}}{2l},\ \ \nu^2=\frac{\varepsilon\lambda\mp\sqrt{\lambda^2-4l(n-l)c^2}}{2(n-l)},
\end{equation*}
\begin{equation}\label{H2}
H^2=\frac{1}{2n^2}\left[n\varepsilon\lambda-4l(n-l)c\varepsilon\pm(2l-n)\sqrt{\lambda^2-4l(n-l)c^2}\right].
\end{equation}

Specially, for the case $\varepsilon\lambda=n|c|$, 
 we check the values in the above equations:
\begin{itemize}
  \item when $c\varepsilon>0$,
  \begin{itemize}
  \item[$\bullet$] $\mu^2=\nu^2=c\varepsilon $ and $H^2=\frac{ (2l-n)^2c\varepsilon}{n^2}
$;
  \item[$\bullet$] $\mu^2=\frac{n-l}{l}c\varepsilon$, $\nu^2=\frac{l}{n-l}c\varepsilon
$ and $
H^2=0
$, a contradiction.
  \end{itemize}
  \item when $c\varepsilon<0$, we have
 \begin{itemize}
  \item[$\bullet$] $\mu^2=-\frac{n-l}{l}c\varepsilon$, $\nu^2=-\frac{l}{n-l}c\varepsilon
$ and $
H^2=-\frac{4l(n-l)}{n^2}c\varepsilon
$;
  \item[$\bullet$] $\mu^2=\nu^2=H^2=-c\varepsilon$, which contradicts to $\mu\neq \nu$.
\end{itemize}
\end{itemize}

By the way, we verify that the right-hand side of \eqref{H2} for the case $\varepsilon\lambda\neq n|c|$ is larger than $0$ and less than $\frac{\varepsilon\lambda}{n}$.
In other words, the values of $H^2$ are all being in the range given by Theorem 3.3.
Indeed, by using $n\varepsilon\lambda\geq2n\sqrt{l(n-l)}|c|\geq4l(n-l)|c|$ and
$$
\frac{[n\varepsilon\lambda\pm4l(n-l)c\varepsilon]^2}{(n-2l)^2[\lambda^2-4l(n-l)c^2]}=1+\frac{4l(n-l)(\lambda\pm nc)^2}{(n-2l)^2[\lambda^2-4l(n-l)c^2]}> 1,
$$
we easily derive
$
 n\varepsilon \lambda\pm4l(n-l)c\varepsilon>\pm(2l-n)\sqrt{\lambda^2-4l(4-l)c^2},
$
which leads to
$$
0<\frac{1}{2n^2}[n\varepsilon\lambda-4l(n-l)c\varepsilon\pm(2l-n)\sqrt{\lambda^2-4l(4-l)c^2}]<\frac{\varepsilon\lambda}{n}.
$$

(i\!i)\ \ When $\mu$ and $\nu$ are imaginary, we suppose $\mu=\gamma+\tau\sqrt{-1}$ and $\nu=\gamma-\tau\sqrt{-1}$.
\eqref{trAlmun-lnu}, \eqref{varepsilonlambdalmun-lnu} and $c+\varepsilon\mu\nu=0$ be rewritten as
\begin{equation*}\label{houbu-eq1-theorem 5.16}
\gamma=\varepsilon H,\  \varepsilon \lambda=n(\gamma^2-\tau^2),\
-c\varepsilon=\gamma^2+\tau^2>0.
\end{equation*}
Then,
$
-nc\varepsilon-\varepsilon \lambda=2n\tau^2>0\ \text{and}\
\varepsilon \lambda-nc\varepsilon=2n\gamma^2>0
$ follows.
Therefore, $|\varepsilon\lambda|<-nc\varepsilon$ and
$
H^2=\gamma^2=\frac{\varepsilon\lambda-nc\varepsilon}{2n}.
$
$\hfill\square$


\section{Classification of non-minimal Lorentzian PMCV hypersurfaces in $\mathbb{H}^{n+1}_1$ and $\mathbb{S}^{n+1}_1$}

In this section, applying the results in Section 3 and the classification theory of isoparametric hypersurfaces,
we classify non-minimal Lorentzian PMCV hypersurfaces of $\mathbb{H}^{n+1}_1(-1)$ and $\mathbb{S}^{n+1}_1(1)\ (n\geq 3)$
with at most two distinct principal curvatures. The classification results for $n=2$ have been provided in \cite{Du 2018}.

According to \cite{Magid 1985}, for Lorentzian hypersurface $M^n_1$ in $N_1^{n+1}(c)$, the shape operator $A$ and the metric matrix $G$
have the following four possible forms:

\begin{footnotesize}
\begin{equation*}
\begin{array}{|c|c|c|c|}
 \hline \! \textbf{Form I}\! & \! \textbf{Form I\!I}\! & \! \textbf{Form I\!I\!I}\! &\! \textbf{Form I\!V}\\
 \hline
 \!A=\!\left(\!
  \begin{array}{ccc}
    \! \lambda_1 \! & \! & \!\\
    \! \! & \!\ddots\! &\! \\
   \! & \! &\! \lambda_{n}\!
  \end{array}\!
\right) \! &\!
 A=\!\left(\!
  \begin{array}{ccc}
    \! \lambda_1 \! & \! & \!\\
    \! 1\! & \!\lambda_1\! &\! \\
   \! & \! &\! D_{n-2}\!
  \end{array}\!
\right)\! &\!
 A=\!\left(\!
  \begin{array}{cccc}
    \! \lambda_1 \!& \! & \! &\!\\
    \! 1 \!& \!\lambda_1\! & \! &\!\\
    \! &\! 1\! &\! \lambda_1\! &\! \\
    \! & \!  &  \!  & \! D_{n-3}\!
  \end{array}\!
\right)\! & \! A=\!\left(\!
  \begin{array}{ccc}
  \! D_{n-2}\! & \! & \!\\
  \! & \!\gamma\! &\! \tau\!\\
  \! &\! -\tau\! & \!\gamma \!
  \end{array}\!
\right)\!
\\
 \!G\!=\!\left(\!\begin{array}{cc}
   \!I_{n-1}\! & \! \\
   \! & \!-1 \!
  \end{array}\!
 \right)\! & \!
   G\!=\!\left(\!
  \begin{array}{ccc}
  \! 0 \!&\! 1\! & \!\\
 \!  1 \!&\! 0\! & \!\\
   \!& \!& \!I_{n-2} \!
  \end{array}\!
\right)\! &\!
   G\!=\!\left(\!
  \begin{array}{cccc}
    \! 0\!& \!0\! & \!1 \!& \!\\
    \!0\! &\! 1 \!&\!0 \! &\! \\
   \! 1 \!&\!0\!  &\!0\! &\!\\
    \!& \! & \!  & \!I_{n-3} \!
  \end{array}\!
\right)\! &\!
G\!=\!\left(\!
  \begin{array}{ccc}
    \!I_{n-2}\! &\! &\!   \\
      \!  & \!   1\! &\! 0\!\\
      \! & \!0 \!&-1\!
  \end{array}\!
\right)\!\\
\hline
\end{array},
\end{equation*}
\end{footnotesize}where $D_{n-2}\!=\!\text{diag}\{\lambda_2, \cdots, \lambda_{n-1}\}$,
$D_{n-3}\!=\!\text{diag}\{\lambda_2, \cdots, \lambda_{n-2}\}$ and $I$ the identity matrix.
For simplicity, a hypersurface whose shape operator has the form I (I\!I, I\!I\!I or I\!V)
is called the type I (I\!I, I\!I\!I or I\!V) hypersurface.

Observe the form I\!V, we find $A$ has at least three distinct eigenvalues when $n\geq 3$.
For the form I, the classification result of non-minimal PMCV hypersurfaces
with at most two distinct principal curvatures has been obtained in \cite{Du 2016}.
So, we only need to discuss the type I\!I and I\!I\!I hypersurfaces.

From Theorem 3.1, we know that PMCV hypersurfaces of $\mathbb{H}^{n+1}_1(-1)$ and $\mathbb{S}^{n+1}_1(1)$
with at most two distinct principal curvatures are isoparametric.
In turn, we can show that the following type I\!I and I\!I\!I isoparametric hypersurfaces
in $\mathbb{H}^{n+1}_1(1)$ and $\mathbb{S}^{n+1}_1(1)$ (cf. \cite{Li 2018, xiao1999}) are also PMCV hypersurfaces
with at most two distinct principal curvatures.

\vskip.2cm
{\bf Example 4.1}\ (Type I\!I hypersurfaces in $\mathbb{H}^{n+1}_1(-1)$, cf. \cite{xiao1999})

Let $\{E_1(t),\ E_2(t),\ \cdots,\ E_{n+2}(t)\}$ be a single parameter frame of $\mathbb{E}^{n+2}_2$ defined on the interval $I$ 
with all scalar products zero except $\langle E_1, E_2\rangle=\langle E_i, E_i\rangle=1\ (3\leq i\leq n+1)$ and $\langle E_{n+2}, E_{n+2}\rangle=-1$,
 and satisfies
$$
\begin{aligned}
&E_1=E_{n+2}^{'},\ E_{n+1}^{'}=\mu E_1+BE_2,\\
&E_{\alpha}^{'}\in \text{span}\{E_2, E_{p+1},\cdots, E_n\},\ p+1\leq\alpha\leq n,
\end{aligned}
$$
where $2\leq p\leq n$, $\mu$ is a nonzero real constant and $B$ is a nonzero function of variable $t$.
The parametrized hypersurfaces in $\mathbb{H}^{n+1}_1(-1)$ given by
\begin{small}
\begin{equation*}
\begin{aligned}
x(t, y)
=&\left(1+\frac{1}{2}\sum_{i=3}^ny_i^2\right)E_{n+2}
-\frac{1}{2}\sum_{i=3}^ny_i^2 E_{n+1}
+\sum_{j=2}^ny_jE_j,\ \text{when}\ \mu^2=1;\\
x(t, y)
=&\text{sgn}(1-\mu^2)\left[\frac{1}{(\mu^2-1)^2}
-\sum_{i=3}^p\frac{y_i^2}{\mu^2-1}\right]^{\frac{1}{2}}
(E_{n+2}-\mu E_{n+1})+\sum_{j=2}^ny_jE_j\\
&+\text{sgn}(\mu(1-\mu^2))\left[\frac{\mu^2}{(\mu^2-1)^2}
+\sum_{\alpha=p+1}^n\frac{y_{\alpha}^2}{\mu^2-1}\right]^{\frac{1}{2}}
(E_{n+1}-\mu E_{n+2})
,\ \text{when}\ \mu^2\neq1,
\end{aligned}
\end{equation*}
\end{small}for $t\in I$ and $y=(y_2,\cdots, y_n)\in \mathbb{E}^{n-1}$, have at most two distinct principal curvatures $\mu$ (with multiply $p$) and $\frac{1}{\mu}$,
and the shape operator $A$ with the form
{\small $$
\left(
\begin{array}{cccc}
\mu  &  & &\\
1 & \mu & &\\
 & & \mu I_{p-2} &\\
  & &  & \frac{1}{\mu}I_{n-p}
\end{array}
\right).
$$}
It follows that $H=\frac{1}{n}\text{tr}A=\frac{p\mu}{n}+\frac{n-p}{n\mu}$ is a real constant
and $\text{tr}A^2=p\mu^2+(n-p)\frac{1}{\mu^2}$. Then, \eqref{proper-eq1} and \eqref{proper-eq2}
hold for $\lambda=p\mu^2+(n-p)\frac{1}{\mu^2}$, i.e. the type I\!I parametrized Lorentzian hypersurfaces
in $\mathbb{H}^{n+1}_1(-1)$ have proper mean curvature vector field.

\vskip.2cm
{\bf Example 4.2}\ (Type I\!I\!I hypersurfaces in $\mathbb{H}^{n+1}_1(-1)$, cf. \cite{xiao1999})

Let $\{E_1(t),\ E_2(t),\ \cdots,\ E_{n+2}(t)\}$ be a single parameter pseudo-orthonormal frame of $\mathbb{E}^{n+2}_2$ defined on the interval $I$ with all scalar products zero except $\langle E_1, E_3\rangle=\langle E_2, E_2\rangle=\langle E_i, E_i\rangle=1\ (3\leq i\leq n+1)$ and $\langle E_{n+1}, E_{n+1}\rangle=-1$,
and satisfies
$$
\begin{aligned}
E_1=E_{n+2}^{'},\
E_{\alpha}^{'}\in \text{span}\{E_3, E_{p+1},\cdots, E_n\},\ p+1\leq\alpha\leq n\ \text{and}\ E_{n+1}^{'}=\mu E_1+BE_2,
\end{aligned}
$$
where $3\leq p\leq n$, $\mu$ is a nonzero real constant and $B$ is a nonzero function of variable $t$.
The parametrized hypersurfaces in $\mathbb{H}^{n+1}_1(-1)$ given by
\begin{small}
\begin{align*}
x(t, y)
=&\left(1+\frac{1}{2}y_2^2+\frac{1}{2}\sum_{i=4}^ny_i^2\right)E_{n+2}
-\frac{1}{2}(y_2^2+\sum_{i=4}^ny_i^2)E_{n+1}
+\sum_{j=2}^ny_jE_j,\ \text{when}\ \mu^2=1;\\
x(t, y)
=&\text{sgn}(1-\mu^2)\left[\frac{1}{(\mu^2-1)^2}
-\frac{y_2^2}{\mu^2-1}-\sum_{i=4}^p\frac{y_i^2}{\mu^2-1}\right]^{\frac{1}{2}}
(E_{n+2}-\mu E_{n+1})+\sum_{j=2}^ny_jE_j\\
&+\text{sgn}(\mu(1-\mu^2))\left[\frac{\mu^2}{(\mu^2-1)^2}
+\sum_{\alpha=p+1}^n\frac{y_{\alpha}^2}{\mu^2-1}\right]^{\frac{1}{2}}
(E_{n+1}-\mu E_{n+2}),\ \text{when}\ \mu^2\neq1,
\end{align*}
\end{small}for $t\in I$ and $y=(y_2,\cdots, y_n)\in \mathbb{E}^{n-1}$, have at most two distinct principal curvatures, and the shape operator takes the form
{\small $$
\left(
\begin{array}{ccccc}
\mu  &  & & &\\
1 & \mu & & &\\
& 1 & \mu & & \\
& & & \mu I_{p-3} &\\
 & & &  & \frac{1}{\mu}I_{n-p}
\end{array}
\right).
$$}
After following the process of Example 4.1, we can also verify that the non-minimal type I\!I\!I parametrized hypersurfaces
in $\mathbb{H}^{n+1}_1(-1)$ are PMCV hypersurfaces,  with the PMCV constant $\lambda=p\mu^2+(n-p)\frac{1}{\mu^2}$.

\vskip.2cm
{\bf Example 4.3}\ (Type I\!I hypersurfaces in $\mathbb{S}^{n+1}_1$(1), cf. \cite{Li 2018})

Let $\{E_1(t),\ E_2(t),\ \cdots,\ E_{n+2}(t)\}$ be a single parameter frame of $\mathbb{E}^{n+2}_1$ defined on the interval $I$ with all scalar products zero except $\langle E_1, E_2\rangle=\langle E_i, E_i\rangle=1\ (3\leq i\leq n+2)$, and satisfies
\begin{equation*}\label{1.5}
\begin{cases}
E_1^{'}=\sum_{i=3}^{n+2} B_{i-2}E_i,\ E_2^{'}=-E_3,\ E_3^{'}=E_1+B_1E_2,\\
E_4^{'}=B_2E_2,\ \cdots,\ E_{n+1}^{'}=B_{n-1}E_2,\ E_{n+2}^{'}=B_nE_2,
\end{cases}
\end{equation*}
where $B_1, B_2, \cdots, B_n$ are smooth functions on $I$ and $\sum_{i=p}^nB_i^2>0$.
For some $2\leq p\leq n$ and $\theta\in[0, 2\pi]$, $\theta\neq\frac{k\pi}{4}\ (k=1, 3, 5, 7)$, the parametrized hypersurface in $\mathbb{S}^{n+1}_1(1)$ given by
$$
x(t,v,y,z)=\cos(\theta-\frac{\pi}{4})\left(vE_2+\sum_{i=3}^{p+1}y_iE_i\right)-\sin(\theta-\frac{\pi}{4})\sum_{\alpha=p+2}^{n+2}z_{\alpha}E_{\alpha},
$$
for $t\in I$, $v\in \mathbb{R}$, $y=(y_3, \cdots, y_{p+1})\in \mathbb{S}^{p-2}_{+}\subset\mathbb{E}^{p-1}$ and $z=(z_{p+2}, \cdots, z_{n+2})\in \mathbb{S}^{n-p}\subset\mathbb{E}^{n-p+1}$, has principal curvatures $\cot(\theta+\frac{\pi}{4})$ (with multiply $p$) and $-\tan(\theta+\frac{\pi}{4})$, and the shape operator with the form
{\small $$
\left(
\begin{array}{cccc}
\cot(\theta+\frac{\pi}{4})  &  &  &\\
1 & \cot(\theta+\frac{\pi}{4}) &  &\\
&  & \cot(\theta+\frac{\pi}{4})I_{p-2} & \\
&  & & -\tan(\theta+\frac{\pi}{4})I_{n-p} \\
\end{array}
\right),\ 2\leq p\leq n.
$$}It follows from $nH=\text{tr}A$ that the mean curvature $H$ is a constant.
Combining $\text{tr}A^2=p\cot^2(\theta+\frac{\pi}{4})+(n-p)\tan^2(\theta+\frac{\pi}{4})$,
it is easy to check that \eqref{proper-eq1} and \eqref{proper-eq2} hold for
$\lambda=p\cot^2(\theta+\frac{\pi}{4})+(n-p)\tan^2(\theta+\frac{\pi}{4})$.
Hence, the non-minimal type I\!I parametrized hypersurface in $\mathbb{S}^{n+1}_1(1)$ has proper mean curvature vector field.

\vskip.2cm
{\bf Example 4.4}\ (Type I\!I\!I hypersurfaces in $\mathbb{S}^{n+1}_1(1)$, cf. \cite{Li 2018})

Let $\{E_1, E_2,\cdots, E_{n+1}\}$ be a parameter frame of $\mathbb{E}^{n+2}_1$ defined on the interval $I$ with all scalar products zero except $\langle E_1, E_3\rangle=\langle E_2, E_2\rangle=\langle E_i, E_i\rangle=1\ (3\leq i\leq n+2)$, and satisfies
\begin{equation*}\label{4.25}
\begin{cases}
E_1^{'}=B_{1}E_2+\sum_{i=4}^{n+2}B_{i-1}E_i,\ E_3^{'}=E_2,\ E_2^{'}=-E_1+B_1E_3,\\
E_r^{'}=B_{r-1}E_3+\sum_{\alpha=p+2}^{n+2}C_{r\alpha}E_{\alpha},\ r=4, \cdots, p+1,\\
E_{\beta}^{'}=B_{\beta-1}E_3-\sum_{s=4}^{p+1}C_{s\beta}E_{s},\ \beta=p+2, \cdots, n+2,\\
\end{cases}
\end{equation*}
where $3\leq p\leq n$, $B_i\ (i=1, 3, \cdots, n), C_{r\alpha}$ ($r=4, \cdots, p+1$, $\alpha=p+2,\ \cdots, n+2$)
are smooth functions on $I$ and the matrix $(C_{r\alpha})\neq0$.
For some $\theta\in[0, 2\pi]$, $\theta\neq\frac{k\pi}{4}\ (k=1, 3, 5, 7)$,
the parametrized hypersurface in $\mathbb{S}^{n+1}_1(1)$ given by
$$
\begin{aligned}
x(t,v,y,z)
=&\cos(\theta-\frac{\pi}{4})\left(vE_3+y_3E_2+\sum_{r=4}^{p+1}y_rE_r\right)-\sin(\theta-\frac{\pi}{4})\sum_{\alpha=p+2}^{n+2}z_{\alpha}E_{\alpha}\\
&-\sqrt{2}\frac{\Psi}{y_3}\sin\theta E_3,\ \Psi=\sum_{r=4}^{p+1}\sum_{\alpha=p+2}^{n+2}y_rz_{\alpha}C_{r\alpha},
\end{aligned}
$$
for $t\in I$, $v\in \mathbb{R}$, $y=(y_3,\cdots, y_{p+1})\in \mathbb{S}^{p-2}\subset\mathbb{E}^{p-1}$ and $z=(z_{p+2}, \cdots, z_{n+2})\in \mathbb{S}^{n-p}\subset\mathbb{E}^{n-p+1}$, has the shape operator of the form
{\small $$
\left(
\begin{array}{ccccc}
\cot(\theta+\frac{\pi}{4})  &  &  & &\\
1 & \cot(\theta+\frac{\pi}{4}) &  & &\\
& 1 & \cot(\theta+\frac{\pi}{4}) & &\\
&  & & \cot(\theta+\frac{\pi}{4})I_{p-3} &\\
& & & & -\tan(\theta+\frac{\pi}{4})I_{n-p} \\
\end{array}
\right).
$$}Similar to Example 4.3, we conclude that the non-minimal type I\!I\!I parametrized hypersurface in $\mathbb{S}^{n+1}_1(1)$
is a PMCV one with the PMCV constant $\lambda=p\cot^2(\theta+\frac{\pi}{4})+(n-p)\tan^2(\theta+\frac{\pi}{4})$.

\vskip.2cm
In the following, we give classification results of non-minimal type I\!I and I\!I\!I PMCV hypersurfaces with at most two distinct principal curvatures. In view of the forms I\!I and I\!I\!I, we know the multiply $l$ of some principal curvature is equal or greater than two or three, respectively. Hence, we can suppose $2\leq l\leq n$ for type I\!I and $3\leq l\leq n$ for type I\!I\!I.

\vskip.2cm
{\bf Theorem 4.5}\quad \emph{Let $M^n_1$ $(n\geq 3)$ be a non-minimal type I\!I (or I\!I\!I)
PMCV hypersurface of $\mathbb{H}^{n+1}_1(-1)$ with the PMCV constant $\lambda$.
Suppose that $M^n_1$ has at most two distinct principal curvatures
with multiplicities $l$ and $n-l$, $2\leq l\leq n$ (or $3\leq l\leq n$),
then $\lambda>0$ for $l=n$, and $\lambda\geq 2\sqrt{l(n-l)}$ for $l<n$.
Furthermore, $M^n_1$ is locally congruent to that of Example 4.1 (or Example 4.2)
with $p=n$, $\mu^2=\frac{\lambda}{n}$ if $l=n$, or $p=l$, $\mu^2=\frac{1}{2l}[\lambda\pm\sqrt{\lambda^2-4l(n-l)}]$ if $l>n$.}

\vskip.2cm
{\bf Proof}\quad We only treat type I\!I hypersurfaces, the other ones are similar. Denote by $\tilde{\mu}$ and $\tilde{\nu}$ the principal curvatures of $M^n_1$ with multiply $l$ and $n-l$.

When $l=n$, i.e. $M^n_1$ has a simple principal curvature, we obtain from Theorem 3.3 that $\lambda>0$ and $H^2=\frac{\lambda}{n}$.
Combining $n\tilde{\mu}=\text{tr}A=n\varepsilon H$, we get $\tilde{\mu}^2=\frac{\lambda}{n}$.
Thus, $M^n_1$ is a totally umbilical isoparametric hypersurface with the square of the principal curvature being equal to $\frac{\lambda}{n}$.
According to \cite{xiao1999}, it is locally congruent to that of Example 4.1 with $p=n$, $\mu^2=\frac{\lambda}{n}$.

When $2\leq l<n$, we easily find from the form I\!I that the geometric multiplicity of $\tilde{\mu}$ is equal to the algebraic multiplicity $n-l$. Applying Theorem 3.5, we derive that $\lambda\geq 2\sqrt{l(n-l)}$ and $\tilde{\mu}^2=\frac{1}{2l}[\lambda\pm\sqrt{\lambda^2-4l(n-l)}]$. So, $M^n_1$ is the isoparametric hypersurface with two distinct principal curvatures (of multiplicities $l$ and $n-l$), and the square of one principal curvature (with multiply $l$) being equal to $\frac{1}{2l}[\lambda\pm\sqrt{\lambda^2-4l(n-l)}]$. The result follows from \cite[Theorem 4.2]{xiao1999}.
\hfill$\square$

\vskip.2cm
Similar to the proof of Theorem 4.5, we obtain the following classification result of non-minimal type I\!I and I\!I\!I PMCV hypersurfaces in $\mathbb{S}^{n+1}_1(1)$.

\vskip.2cm
{\bf Theorem 4.6}\quad \emph{Let $M^n_1$ $(n\geq 3)$ be a non-minimal type I\!I (or I\!I\!I) PMCV hypersurface of $\mathbb{S}^{n+1}_1(1)$ with the PMCV constant $\lambda$.
Suppose that $M^n_1$ has at most two distinct principal curvatures with multiplicities $l$ and $n-l$, $2\leq l\leq n$ (or $3\leq l\leq n$),
then $\lambda>0$ for $l=n$, and $\lambda\geq 2\sqrt{l(n-l)}$ for $l<n$.
Furthermore, $M^n_1$ is locally congruent to that of Example 4.3 (or Example 4.4) with $p=n$, $\cot^2(\theta+\frac{\pi}{4})=\frac{\lambda}{n}$ if $l=n$,
or $p=l$, $\cot^2(\theta+\frac{\pi}{4})=\frac{1}{2l}[\lambda\pm\sqrt{\lambda^2-4l(n-l)}]$ if $l>n$.}


%

\end{document}